\parindent=0pt              
\baselineskip=20pt
\magnification=\magstep1   

\centerline{\bf Small Gamma Products with Simple Values}
\medskip
\centerline{Albert Nijenhuis}
\medskip

{\bf Introduction.} Central as the gamma function is, it is surprising that its only known simple specific values are $\Gamma(m)$ and $\Gamma(m+{1\over2})$ for integral $m$. There are, however, numerous formulas that relate specific values of $\Gamma$ to other functions, such as elliptic or hypergeometric functions - or to other specific values of $\Gamma$ itself. Among the latter there are products  that have very simple values. An example is 
$$\Gamma\left({1\over14}\right)\Gamma\left({9\over14}\right)\Gamma\left({11\over14}\right)=4 \pi^{3/2},\leqno(1)$$
which recently occurred  in the Problems section of the ``Monthly" [1], see also [2]. There is also, of course, the classical multiplier formula
$$\prod_{k=1}^{m-1}\Gamma\left({k\over m}\right) = (2\pi)^{(m-1)/2} m^{-1/2},\leqno(2)$$
but it would be more interesting to find simple values for products of fewer factors.
In this note we do just that, for a large class of products, and with little computational effort. 

Consider any odd integer $n > 1$ and the  set $\Phi(2n)$  of numbers in $[0,2n]$ that are relatively prime to $2n$.  (Its cardinality is $\varphi(2n)$,  Euler's totient function.) $\Phi(2n)$  is a group with respect to  multiplication modulo $2n$. Let $\nu(n)$ be the order of the subgroup generated by $n+2$, and let $A$ be this subgroup or any one of its cosets. Let $b(A)$ count the  $x\in A$ that are larger than $n$. Our main result is
\medskip
{\bf Theorem.}
$$\prod_{x\in A} \Gamma\left({x\over 2n}\right) = 2^{b(A)}\pi^{\nu(n)/2}.\leqno(3)$$
\medskip
A few  cases of this formula are known, but even for small $n$ there are some gaps in published lists.
\bigskip
{\bf The Formula.} 
 Throughout this paper, $n > 1$ will denote a ``fixed" odd integer, and let $\Phi(n)$ be the set (group) of all integers in the interval $[0,n]$  relatively prime to $n$. Then $\phi(n) = \phi(2n)$, and the map $\alpha: \Phi(n) \to \Phi(2n)$ given by
$$ \alpha: \Phi(n) \to \Phi(2n),\quad  \alpha(y) = y \quad\hbox{ if}\quad y \quad\hbox{ odd,  else}\quad \alpha(y)= y+n\leqno(4)$$
is a group isomorphism. (The proof is a simple exercise, distinguishing 3 cases.) The inverse is
$$\alpha^{-1}: \Phi(2n) \to \Phi(n), \quad\alpha^{-1}(x)=x \quad\hbox{if}\quad x<n,\quad\hbox{else}\quad \alpha^{-1}(x)=x-n.$$

We also need a map $\beta$ which ``halves" the elements of $\Phi(n)$  (modulo $n$):
$$\beta: \Phi(n) \to \Phi(n),\quad \beta(y) = y/2 \quad\hbox{ if}\quad  y \quad\hbox{ even,  else}\quad \beta(y)=(y+n)/2.\leqno(5)$$
The doubling formula for $\Gamma$ is needed in the following form
$$ \Gamma(t) = c_t \Gamma(2t)\big/\Gamma\left(t+{1\over2}\right),\leqno(7)$$
 where $c_t=(2\sqrt\pi)2^{-2t}$.
 In (7) set $t=x/2n$, where $x\in\Phi(2n)$,
$$ \Gamma\left({x\over 2n}\right)=c_{x/2n}\Gamma\left({x\over n}\right)/\Gamma\left({x+n\over 2n}\right).\leqno(8)$$
This equation is of the form
$$\Gamma\left({x\over 2n}\right)=\varepsilon(x)c_{x/2n}\Gamma\left({y\over n}\right)/\Gamma\left({z\over n}\right)\leqno(9)$$
with $y=x,z=(x+n)/2 \in \Phi(n)$ and $\varepsilon(x)=1$ when $x < n$. When $x>n$, we apply the reduction formula $\Gamma(t)=(t-1)\Gamma(t-1)$ to both $\Gamma$'s on the right in (8), and cancel factors $(x-n)/n$. The result is
$$\Gamma\left({x\over 2n}\right)=2c_{x/2n}\Gamma\left({x-n\over n}\right)\bigg/\Gamma\left({x-n\over 2n}\right),$$
which is of the form (9) with  $y=x-n, z=(x-n)/2\in \Phi(n)$ and $\varepsilon(x)=2$. 
\medskip
{\bf Lemma.} Let $n$ be an odd integer, $n > 1$
, and $x\in \Phi(2n)$. Then (9) holds, where $y=\alpha^{-1}(x), z=\beta(y) \in \Phi(n)$. Further, $x-n = 2y - 2z$ and $y\equiv 2z \,\,({\rm mod\,\,} n)$.
\medskip
Proof. Distinguish the two cases $x<n$ and $x>n$. (Note that  $x\in\Phi(2n)$ is odd.)

If $x<n$, then $\alpha^{-1}(x) = x = y$ (odd), so $\beta(y)=(y+n)/2=(x+n)/2=z$.

 Also, $2y-2z=2x-(x+n)=x-n$ and $2z=y+n\equiv y$ (mod $n$).

If $x > n$, then $\alpha^{-1}(x)=x-n$ (even), so $\beta(y)=y/2=(x-n)/2=z$

Also, $2y-2z=2(x-n)-(x-n)=x-n$, and $2z=y$.

{\bf Proof of (3).} The members of $\Phi(n)$ are taken as vertices of a directed labeled graph. The edges are the pairs $(y,z)=(y,\beta(y))$; such an edge is labeled $x = \alpha(y)\in \Phi(2n)$. The vertices have in- and outdegree 1, so the connected components are cycles. In fact, the component of 1 is the cyclic group generated by 2; denote its order by $\nu(n)$. The other components are the cosets. Let $B$ be any one of these cycles. Form the product $P = \prod\Gamma(x/2n)$ of the left sides of (9), where the product extends over $x\in \alpha B$, i.e.,the labels of the edges of $B$. Similarly, take the product of the right sides of (9), and note that the product telescopes as all the  $\Gamma$'s cancel, leaving only
$$ P = \prod_{\alpha B}\varepsilon(x)c_{-x/2n}=(2\sqrt\pi)^{\nu(n)}2^{-2(\Sigma x)/2n}\prod_{\alpha B}\varepsilon(x).$$
Similarly, we have the telescoping sum $\sum(x-n) = \sum(2y-2z)=0$, so $\sum x = n\nu(n)$. Therefore,
$$ P=(2\sqrt\pi)^{\nu(n)}2^{-2n\nu(n)/2n}\prod_{x\in\alpha B}\varepsilon(x) = \pi^{\nu(n)/2}2^{b(\alpha B)},$$
where $b(\alpha B)$ is the number of $x \in \alpha B$ that are bigger than $n$. In summary,
$$\prod_{x\in\alpha B}\Gamma\left({x\over 2n}\right) = 2^{b(\alpha B)} \pi^{\nu(n)/2}.$$
 Since $\alpha$ is an isomorphism,  $A=\alpha B$  is (a coset of) the subgroup of $\Phi(2n)$ generated by $\alpha(2) = n+2$. That yields (3).
\bigskip
{\bf Corollaries.}

{\bf 1.} In Zucker [2] we find,  for $n=2^m-1$, $m>1$, that
$$\Gamma\left({1\over 2n}\right)\prod_{k=1}^{m-1}\Gamma\left({2^k+n\over2n}\right)=2^{m-1}\pi^{m/2}.$$
Proof: This is a special case of (3). The left side equals
$$\prod_{k=0}^{m-1}\Gamma\left({\alpha(2^k)\over 2n}\right),$$
while $\beta(1)=2^{m-1}, \beta(2^k) = 2^{k-1} \, (k>0)$. All numerators except for one are $>n$.

{\bf 2.} Complementation. Let $A$ be as in (3), and $A^* = \{2n-x|x\in A\}$, then (3) holds with $A^*$ replacing $A$, and $b(A^*)=\nu(n)-b(A)$.

Proof: Verify that if $x'\equiv x(n+2) \pmod{2n}$ then $2n-x'\equiv (2n-x)(n+2) \pmod{2n}$. 

{\bf 3.} Take the product in (3) over all of $\Phi(2n)$; that is, multiply both sides of (3) as $A$ ranges over the subgroup and its cosets. There are $\varphi(n)/\nu(n)$ choices of $A$. Each $x\in \Phi(2n)$ greater than $n$ will occur exactly once. 
$$ \prod_{x\in\Phi(2n)} \Gamma\left({x\over 2n}\right) = (2\pi)^{\varphi(n)/2}.$$

{\bf 4.}  Some numerical examples. We determined, for which odd $n<100$, the number of  sets  $A$ exceeds 2. There are 9 such values. Among these, only $n=43$ has 3 $A$'s, and all of these
 are self-complementary [2].  In all other cases the number is even, and can be as big as 8. At the other end, there are 16 values of $n$ for which $\nu(n)=\phi(n)$; that can only happen when $n$ is a prime or a prime power. 

We list the six sets $A$ for $n=31$ because only two of them are usually mentioned [2]. 
\medskip
$$\eqalign{
n=31: &(1,33,35,39,47), (3,17,37,43,55), (5,9,41,49,51),\cr
 &(7,19,25,45,59), (11,13,21,53,57), (15,23,27,29,61).}$$
The first one, written out in full, is
$$  \Gamma\left({1\over62}\right)\Gamma\left({33\over62}\right)\Gamma\left({35\over62}\right)\Gamma\left({39\over62}\right)\Gamma\left({47\over62}\right)=2^4\pi^{5/2}.$$
Here, $\nu=5$, the length of the product, and $b=4$, the number of numerators bigger than 31. Each numerator, multiplied by 33 (mod 62) yields the next one, in circular order.
\medskip
In a personal note,  H. Chen states the problem of finding minimum sizes of gamma products that have simple values. This paper may be a step in that direction, but any definitive answer will depend on a suitable definition of ``simple value".

\medskip
{\bf References}

[1] Glasser, M. L.,  Problem 11426, Amer. Math. Monthly, {\bf 116}, p. 365, (2009)

[2] Zucker, L.J., Personal notes (1994)
\medskip
Department of Mathematics, University of Pennsylvania

Department of Mathematics, University of Washington

nijenhuisalbert@msn.com

\bye